\documentclass[reviewcopy]{elsart}
\usepackage{amssymb}
\usepackage{amssymb}
\usepackage{amssymb}
\usepackage{amssymb}
\usepackage{amssymb}
\usepackage{amsfonts}
\usepackage{amssymb}
\usepackage{amssymb}
\usepackage{amssymb}
\usepackage{amssymb}

\usepackage{amssymb}
\usepackage{amsfonts}
\usepackage{amsmath}
\usepackage{fancyhdr}

\begin{document}

\begin{frontmatter}



\title{On the Ring  of Integer-valued Quasi-polynomials\thanksref{titlefn}}
\thanks[titlefn]{Projects 10526016
 and 10571033 Supported by National Natural
Science Foundation of China. Project Supported by Development
Program for Outstanding Young Teachers in Harbin Institute of
Technology.}

\author{Nan Li}\ead{amenda860111@gmail.com },
\author{Sheng Chen\corauthref{cor}}
\corauth[cor]{Corresponding author.} \ead{schen@hit.edu.cn}
\address{Department of Mathematics,   Harbin Institute of Technology\\
Harbin, 150001,   P. R. China}

\begin{abstract}
The paper studies some properties of the ring  of integer-valued
quasi-polynomials.  On this ring,   theory of  generalized Euclidean
division and generalized greatest common divisor are presented.
Applications to finite simple continued fraction expansion of
rational numbers and Smith normal form of integral matrices with an
integer parameter are also
given. 
\end{abstract}

\begin{keyword}
ring of integer-valued quasi-polynomials,  generalized Euclidean
division, continued fraction, Smith normal form, integer parameter


\MSC 11A05;13F20;05A15

\end{keyword}
\end{frontmatter}

\label{}


\section{Introduction}
Integer-valued quasi-polynomials occur naturally in graded algebras
(see \cite{Atiyah}) and enumerative combinatorics (see
\cite{stanley}),   such as the Ehrhart quasi-polynomial of a
rational polytope.  In this article,   we study division theory of
the ring $R$ of all integer-valued quasi-polynomials and its
applications.

Obviously  $\mathbb{Z}\varsubsetneq \mathbb{Z}[x]\varsubsetneq R $,
where $\mathbb{Z}[x]$ is a unique factorization domain(UFD) but not
an Euclidean Domain. It turns out that  $R$ is neither a domain nor
a Noetherian ring, but every finitely generated ideal in  $R$ is a
principal ideal(see Corollary \ref{cor1} ). However, we can develop
theory of generalized Euclidean division and greatest common
divisor(GCD) on $R$, which has close relations to Euclidean division
and GCD theory on $\mathbb{Z}$ pointwisely(i.e., through
evaluation).




The organization of the paper is as follows.  In Section 2,   we
will give the definition of integer-valued quasi-polynomial (see
Definition $\ref{2. 0. 002}$) and prove some elementary results.  In
Section 3, we will present generalized Euclidean division on $R$.
The relation between this generalized division and Euclidean
division on $\mathbb{Z}$ is considered in Remark $\ref{2. 0. 01}$.
In Section 4, we will present generalized GCD theory through
successive generalized Euclidean divisions. Then, in Section 5  we
will develop generalized GCD theory through pointwisely GCD over
$\mathbb{Z} $ and  prove the equivalent of the generalized GCD
developed in Section 4 and 5.(see Theorem $\ref{2. 0. 015}$). In
Section 6, we will give some applications of our generalized
Euclidean division and GCD theory.
 We will show that
expanding rational functions $h(x)=\frac{f(x)}{g(x)}$,   where
$f(x), g(x)\in \mathbb{Z}[x]$,   into finite simple continued
fractions for every $x\in \mathbb{Z}$,   the numbers of terms for
$h(x)$ are uniformly bounded (see Theorem $\ref{finite}$).
Applications in Smith normal form
 for integral matrices with an integer parameter are also given( see Theorem $\ref{3. 0. 6}$). \par Most of
 our
proofs are simple but constructive.  We have implemented all of the
related algorithms  in Maple and would like to send the Maple files
to the interested readers upon request.

\section{Ring of integer-valued Quasi Polynomial}

We begin with two simple lemmas.
\begin{lem}\label{2. 0. 000}  Suppose that  $f(x)\in \mathbb{Q}[x]$ and for every $n\in
\mathbb{Z}$,   $f(n)\in \mathbb{Z}$.  Let $a$ be a positive integer
such that $af(x)\in \mathbb{Z}[x]$.  Then for every $i=0,  1,
\cdots, a-1$,   we have
$$g_i(x)=f(ax+i)\in \mathbb{Z}[x]$$
\end{lem}

\begin{pf} Let $f(x)=b_lx^{l}+b_{l-1}x^{l-1}+\cdots +b_1x+b_0$.
Then for $i=0,  1,  \cdots,  a-1$,  using  binomial expansion
theorem, we have
\begin{eqnarray*}\label{2. 0. 010}
  f(ax+i) &=& b_l(ax+i)^{l}+b_{l-1}(ax+i)^{l-1}+\cdots+b_1(ax+i)+b_0 \\
   &=& b_l[(ax)^{l}+\cdots+C_{l}^{1}(ax)i^{l-1}+i^{l}]\\&&+b_{l-1}[(ax)^{l-1}+\cdots+C_{l-1}^{1}(ax)i^{l-2}+i^{l-1}]\\&&+\cdots+b_1(ax+i)+b_0
\end{eqnarray*}
So there exists $h_i(x)\in \mathbb{Q}[x]$ such that
$$f(ax+i)=ah_i(x)+f(i)$$
Since $af(x)\in \mathbb{Z}[x]$,   we get $ah_i(x)\in \mathbb{Z}[x]$.
   For $i=0,  1, \cdots,  a-1$,  $f(i)\in
\mathbb{Z}$ and   thus we have $g_i(x)=f(ax+i)\in \mathbb{Z}[x]$.
\qed
\end{pf}

\begin{lem}\label{2. 0. 001}For a function $f:\mathbb{Z}\to\mathbb{Z}$,   the following
two conditions are equivalent:
\begin{enumerate}
  \item[(1)] there exists  $T\in \mathbb{N}$ and $f_i(x)\in \mathbb{Q}[x]$
  ($i=0,  1,  \cdots,  T-1$) such that for any $m\in \mathbb{Z}$,
  if $n=Tm+i$,   we have $f(n)=f_i(m)$;
  \item[(2)] there exists $T\in \mathbb{N}$ and $f_i(x)\in \mathbb{Z}[x]$($i=0,  1,  \cdots,  T-1$) such that for any $m\in
  \mathbb{Z}$,
  if
  $n=Tm+i$,   we have $f(n)=f_i(m)$.
\end{enumerate}
\end{lem}
\begin{pf}
It suffices to show that $(1)$ implies $ (2)$.  Suppose that there
exists $T_0\in \mathbb{N}$ and $g_i(x)\in \mathbb{Q}[x]$($i=0,  1,
\cdots, T_0-1$) such that for every $n\in \mathbb{Z}$,   when
$n=T_0k+i$, we have $f(n)=g_i(k)$.  Choose an integer $a$ large
enough such that for each $i\in \{0,  1,  \cdots, T_0-1\}$,
$ag_i(x)\in \mathbb{Z}[x]$. Define
$$
  f_{ij}: \mathbb{Z} \rightarrow \mathbb{Z}:\\
  m \longmapsto  g_{i}(am+j)
$$
where $i=0,  1,  \cdots,  T_0-1$ and  $j=0,  1,  \cdots,  a-1$.
Since $ag_i(x)\in \mathbb{Z}[x]$ and $g_i(j)=f(Tj+i)\in \mathbb{Z}$,
by Lemma $\ref{2. 0. 000}$,   for each $j\in\{0,  1,  \cdots,
a-1\}$, we have $f_{ij}(x)\in \mathbb{Z}[x]$.  Now let $T=aT_0$. For
$l=jT_0+i\in \{0,  1,  \cdots,  T-1\}$,   define $f_l(x)=f_{ij}(x)$.
Then $f(x)$ satisfies condition (2). \qed
\end{pf}

\begin{defn}\label{2. 0. 002}We call  function $f(x)$ an \emph{integer-valued
quasi-polynomial}   if it satisfies the equivalent conditions in
Lemma $\ref{2. 0. 001}$.  If $T\in \mathbb{N}$ and $f_i(x)\in
\mathbb{Z}[x]$ ($i=0,  1,  \cdots,  T-1$) satisfy  condition 2 in
Lemma $\ref{2. 0. 001}$,   we call $(T,  \{f_i(x)\}_{i=0}^{T-1})$ a
\emph{representation} of $f(x)$ and write
$$f(x)=(T,  \{f_i(x)\}_{i=0}^{T-1})$$We call $\max\{degree(f_i(x))|
i=0,  1,  \cdots,  T-1\}$ and $T$ the  degree and period of this
representation respectively.
\end{defn}We can see that degree of $f(x)$ is independent from its representations.  In fact,   let
$(T_1,  \{f_i(x)\}_{i=0}^{T_1-1})$ and $(T_2,
\{g_j(x)\}_{j=0}^{T_2-1})$ be two representations of $f(x)$.  For
each $s\in \{0,  1,  \cdots,  \frac{T_1T_2}{d}-1\}$,   if
$s=T_1s_1+i=T_2s_2+j$,   where $s_1,   s_2\in \mathbb{Z} $,   $ 0
\leq i< T_1 $ and $ 0 \leq j< T_2 $, then
$\frac{T_1T_2}{d}x+s=T_1(\frac{T_2}{d}x+s_1)+i=T_2(\frac{T_1}{d}x+s_2)+j$,
and therefore  we have
$$f_i(\frac{T_2}{d}x+s_1)=g_j(\frac{T_1}{d}x+s_2)$$
Hence it is not difficult to check that
$$max\{degree(f_i(x))|
i=0,  1,  \cdots,  T_1-1\}=max\{degree(g_i(x))| i=0,  1,  \cdots,
T_2-1\}$$  The periods are not unique.  However, we have the
following result.
\begin{prop}\label{2. 0. 011} Let $f(x)$ be an integer-valued
quasi-polynomial and $ \Omega $ be the set of positive periods $T$
satisfying the equivalent conditions in Lemma $\ref{2. 0. 001}$.
Define a partial order $\preccurlyeq\,  \,  $ on $ \Omega $ as
follows: $\,  T_1\preccurlyeq T_2$ if and only if  $T_1\mid T_2$.
Then $ \Omega $ has a smallest item,   which will be called the
least positive period of $f(x)$.
\end{prop}
\begin{pf} It is sufficient for us to prove that if $f(x)\neq{0}$ and $ T_1,  T_2\in
\Omega $,   then $d=gcd_{\mathbb{Z}}(T_1, T_2)\in
                \Omega $,   i.e.,  if
$f(x)=(T_1,  \{f_i(x)\}_{i=0}^{T_1-1})=(T_2,
\{g_j(x)\}_{j=0}^{T_2-1})\neq 0$,  then $f(x)$ has a representation
with period $d$.  Suppose that $k$ is the degree of  $f(x)$.  Let $$
f(T_1x+i) =f_i(x)=a_ix^{k}+\cdots
$$
  and $$
  f(T_2x+j) = g_j(x)=b_jx^{k}+\cdots
$$

Now we are going to show $d\in  \Omega $ by induction on $k$. For
$k=0$,   the assertion reduces to $f(n)=f(n+d)$,   for every $ n\in
\mathbb{Z}$.  Note that there exist $u,  v\in \mathbb{Z}$ such that
$T_1u+T_2v=d$.  Since $ T_1,  T_2\in
 \Omega $,   for every $x\in
               \mathbb{Z}$,   we have $f(x+d)=f(x+T_1u+T_2v)=f(x+T_2v)=f(x)$,   as desired.

Suppose the assertion true for $k-1$ and prove it for $k$.  Assume
that $s\in \{0,  1,  \cdots,  \frac{T_1T_2}{d}-1\}$,   and
$s=T_1s_1+i=T_2s_2+j$,   where $s_1,   s_2\in \mathbb{Z} $,   $ 0
\leq i< T_1 $ and $ 0 \leq j< T_2 $.  Then since
$\frac{T_1T_2}{d}x+s=T_1(\frac{T_2}{d}x+s_1)+i=T_2(\frac{T_1}{d}x+s_2)+j$,
we have
$$f_i(\frac{T_2}{d}x+s_1)=g_j(\frac{T_1}{d}x+s_2)$$
Take  the $k$-th order derivative to the above equation,   we get
$$
  k!(\frac{T_2}{d})^{k}a_i =f_i^{(k)}(\frac{T_2}{d}x+s_1) =
 g_j^{(k)}(\frac{T_1}{d}x+s_2) =k!(\frac{T_1}{d})^{k}b_j
$$
The left and right of the identity are all integer constant function
with periods $T_1$ and $T_2$ respectively.  By the assertion for
$k=0$,   when $i\equiv j\,  \,   (mod\,  \,  T)$,   we have
$$(\frac{T_2}{d})^{k}a_i =(\frac{T_2}{d})^{k}a_j= (\frac{T_1}{d})^{k}b_i=(\frac{T_1}{d})^{k}b_j$$
  Define
$$h:\mathbb{Z}  \rightarrow \mathbb{Z}:n \mapsto a_im^{k} $$ where $n=T_1m+i$,
$i=0,  1,  \cdots,   T_1-1$. Let $l(x)=f(x)-h(x)$.  Then $l(x)\in R$
has periods $T_1$,   $T_2$ and
 $k-1$ is its  degree.  By induction hypothesis,   $l(x)\in
R$ and has a period $d$.  Note that $h(x)$ has a period $d$.
Consequently, $f(x)=l(x)+h(x)$ has a period $d$,   as desired.
\qed\end{pf}
\begin{rem}\label{2. 0. 003}
In the remainder of this paper,   we often use the following fact.
Let $f(x),  g(x)$ be two integer-valued quasi-polynomials with
periods,   say,   $T_1,   T_2$ respectively. Then $f(x),  g(x)$ have
a common period  $T=lcm(T_1,   T_2)$.
\end{rem}

\begin{prop}\label{R}The set of all integer-valued quasi-polynomials,   denoted
by $R$,   with pointwisely defined addition and multiplication, is a
commutative ring with identity.
\end{prop}
\begin{pf}
Let $\Gamma$ be the set of functions $f:\mathbb{Z}\to\mathbb{Z}$. It
is obvious that  $1\in R$ and with pointwisely defined addition and
multiplication, $\Gamma$ is a commutative ring with identity $1$.
Therefore,   it is sufficient for us to prove $R$ is a subring of
$\Gamma$. In fact, by remark $\ref{2. 0. 003}$, we know that
subtraction and multiplication are closed in $R$.  The proof is
completed. \qed
\end{pf}

\begin{prop}\label{2. 0. 004}

\begin{enumerate}
  \item[(1)] Let $f(x)\in R$ and $f(x)\neq 0$.  Then f(x) is a zero divisor if and only
  if it has a representation $(T,  \{f_i(x)\}_{i=0}^{T-1})$
  and there exists
  $i\in\{0,  1,  \cdots,  T-1\}$ such that $f_i(x)=0$.
  \item[(2)] Let $f(x)\in R$.  Then f(x) is  invertible  if and only if it has
  a representation $(T,  \{f_i(x)\}_{i=0}^{T-1})$ such that for
  $i\in\{0,  1,  \cdots,  T-1\}$,   $f_i(x)=1$ or $f_i(x)=-1$.
\end{enumerate}
\end{prop}
\begin{pf} (1) Suppose $f(x)g(x)=0$ where  $g(x)\in R- \{0\}$.  Let $T$ be a
common period of $f(x)$ and $g(x)$ such that $f(x)=(T,
\{f_i(x)\}_{i=0}^{T-1})$ and $g(x)=(T,  \{g_i(x)\}_{i=0}^{T-1})$
(see Remark $\ref{2. 0. 003}$).  Since $g(x)\neq 0$,   there exists
$i_0\in \{0,  1,  \cdots,  T-1\}$ with $g_{i_0}(x)\neq 0$ and thus
$f_{i_0}(x)= 0$. \par Conversely,  suppose that there exists $i_0\in
\{0,  1,  \cdots,  T-1\}$ with $f_{i_0}(x)=0$.  Then define $g(x)\in
R$ as follows: if $n=Tm+i_0$,  $g(n)=g_{i_0}(m)=m$,   and else
$g(n)=0$. Obviously,   we have $f(x)g(x)=0$.  Note that $f(x)\neq
0$.  Hence $f(x)$ is a zero divisor in $R$.

 (2) Assume that there exists $g(x)\in R-{0}$ with $f(x)g(x)=1$.
Let $T$ be a common period of $f(x)$ and $g(x)$,   such that
$f(x)=(T,  \{f_i(x)\}_{i=0}^{T-1})$ and $g(x)=(T,
\{g_i(x)\}_{i=0}^{T-1})$.  Since for each $i\in \{0,  1,  \cdots,
T-1\}$,   $f_i(x)g_i(x)=1$ and $f_i(x),  g_i(x)\in \mathbb{Z}[x]$.
It follows that $f_i(x)=g_i(x)=1$ or $-1$.  The reverse part is
trivial. \qed
\end{pf}

\begin{prop}\label{2. 0. 005}The ring $R$ is not Noetherian.
\end{prop}
\begin{pf}
Let$f_1(n)=n$ and $$f_2(n)=\left\{\begin{array}{ll}
            \frac{n}{2},   & \hbox{$n=2m$} \\
            n,   & \hbox{$n=2m+1$}
            \end{array}
        \right. ,
f_3(n)=\left\{
          \begin{array}{ll}
            \frac{n}{4},   & \hbox{$n=4m$} \\
            \frac{n}{2},   & \hbox{$n=4m+2$} \\
            n,   & \hbox{else}
          \end{array}
        \right. ,  \cdots
        f_k(n)=\left\{\begin{array}{ll}
            \frac{f_{k-1}(n)}{2},   & \hbox{$n=2^{k-1}m$} \\
            f_{k-1}(n),   & \hbox{else}
            \end{array}
        \right. \cdots$$
Let
$$g_2(n)=\left\{\begin{array}{ll}
            2,   & \hbox{$n=2m$} \\
            1,   & \hbox{$n=2m+1$}
            \end{array}
        \right. ,   g_3(n)=\left\{\begin{array}{ll}
            2,   & \hbox{$n=4m$} \\
            1,   & \hbox{else}
            \end{array}
        \right. ,   \cdots
        g_k(n)=\left\{\begin{array}{ll}
            2,   & \hbox{$n=2^{k-1}m$} \\
            1,   & \hbox{else}
            \end{array}
        \right. \cdots$$
Note that $f_{k-1}(n)=f_k(n)g_k(n)$.  We have
$(f_{k-1}(n))\varsubsetneq (f_{k}(n))$,  where $k=2,  3, \cdots$.
Thus $R$ does not satisfy ascending chain condition(acc) on ideals
in $R$,  i.e., it is not Noetherian(cf. \cite{Ring}). \qed
\end{pf}

\section{Generalized Euclidean Division}

\begin{defn}\label{2. 0. 2}
               Let $r(x)\in R$ and $r(x)=(T,
 \{r_i(x)\}_{i=0}^{T-1})$. We shall say $r(x)$ is nonnegative and  write $r(x)\succcurlyeq  0$,   if it
satisfies the following equivalent conditions ：
\begin{enumerate}
\item[(1)]
for every $i=0,  1,  \cdots,  T-1$, $r_i(x)=0$ or its leading
coefficient is positive;
          \item[(2)] there exists $C\in \mathbb{Z}$,   such that for every integer $n>C$,  we have $r(n) \geqslant 0$.
          \end{enumerate}
We shall say $r(x)$ is strictly positive and  write $r(x)\succ 0$,
if $r(x)=(T,
 \{r_i(x)\}_{i=0}^{T-1})$ satisfies the following  condition: $(1^{'})$ for every $i=0,  1,  \cdots, T-1$,
the leading coefficient of $r_i(x)$ is positive. {\bf We write $
f(x)\preccurlyeq   g(x)$  if $ g(x)-  f(x)\succcurlyeq 0$}.
\end{defn}
\begin{rem}
Because of the existence of zero divisor, the "order" in $R$ defined
above is not a partial order, i.e., $r(x)\in R$ satisfies both
$r(x)\succcurlyeq 0$ and $r(x)\preccurlyeq 0$ does not imply
$r(x)=0$. More precisely,  $r(x)\succcurlyeq 0$ implies $r(x)$ falls
in one of the following three cases: (a) $r(x)\succ 0$, (b) $r(x)=
0$, (c) $r(x)$ is a nonnegative zero divisor, where case $(c)$ is
impossible if  $r(x)\in  \mathbb{Z}[x]$(by Proposition $\ref{2. 0.
004}$). Thus the following definition is {\bf well defined}.
\end{rem}.
\begin{defn}\label{2. 0. 20}
Let $r(x) \in \mathbb{Z}[x] $,   define a function $| . |$ as
follows: $\mathbb{Z}[x]\rightarrow \mathbb{Z}[x]$:
$$|r(x)|=\left\{
          \begin{array}{ll}
            r(x),   & if\,  \,  \hbox{$r(x)\succ 0$} \\
            -r(x),   & if \,  \,   \hbox{$r(x)\prec 0$} \\
            0,   & if \,  \,  \hbox{$r(x)=0$}
          \end{array}
        \right. $$
\end{defn}

\begin{thm}\label{2. 0. 3}
Let $f(x),  g(x)\in \mathbb{Z}[x]$ and $g(x)\neq 0$.  Then there
exist unique $P(x),  r(x)\in R$ such that
$$f(x)=P(x)g(x)+r(x)\quad\,  where\quad
0\preccurlyeq r(x)\prec| g(x)|$$ In this situation,  we write
$P(x)=quo(f(x), g(x))$ and $r(x)=rem(f(x),g(x))$.
\end{thm}

\begin{pf} \emph{Step 1}: We first prove for the existence of $P(x)$ and $r(x)$.
If  $f(x)=0$, then put $P(x)=0,  r(x)=0$. If $f(x)\neq 0$, then  let
$$ f(x)=a_kx^{k}+a_{k-1}x^{k-1}+\cdots+a_1x+a_0$$
and
$$
g(x)=b_lx^{l}+b_{l-1}x^{l-1}+\cdots+b_1x+b_0
$$
where $a_k\neq 0$ and $b_l\neq 0$.  \par When $k<l$,    define
$P(x), r(x)$ as follows:
\begin{enumerate}
  \item if $f(x)\succ 0$,
then P(x)=0,  r(x)=f(x);
  \item if $f(x)\prec 0$ and $g(x)\prec 0$,
then P(x)=1,  r(x)=f(x)-g(x);
\item if $f(x)\prec 0$ and $g(x)\succ 0$,
then P(x)=-1,  r(x)=f(x)+g(x).
\end{enumerate}
It is easy to check they satisfy the remainder conditions in the
result. When $k\geqslant l$, we claim that there exist $P(n),
r(n)\in R$ with a period $T=b_l^{k-l}$ and prove it by induction on
$k-l$. This assertion is trivial for $k-l=0$. Suppose that it is
true for $k-l=h$ and prove it for $k-l=h+1$.  Now for $i=0, \cdots,
b_l-1$, when $n=b_lm+i$, define
$$h_i(m)=f(b_lm+i)=a_kb_l^{k}m^{k}+h_{i,  k-1}(m)$$
and
$$l_i(m)=g(b_lm+i)=b_l^{l+1}m^{l}+l_{i,  l-1}(m)$$
where $h_{i,  k-1}(x)$ and $ l_{i,  l-1}(x)\in \mathbb{Z}[x]$ have
degrees $k-1$ and  $l-1$ respectively.  Since for every $m\in
\mathbb{Z}$
$$a_kb_l^{k}m^{k}=a_kb_{l}^{h}\,  m^{h+1}\,  l_i(m)-a_kb_l^{h}\,  m^{h+1}\,  l_{i,l-1}(m)$$
we have $$h_i(m)=a_kb_{l}^{h}\,  m^{h+1}\,  l_i(m)+S_{i,  k-1}(m)$$
where$$S_{i,  k-1}(m)=h_{i,  k-1}(m)-a_kb_l^{h}\,  m^{h+1}\,
l_{i，l-1}(m)$$ Note that $S_{i,  k-1}(x)$ and $l_i(x)\in
\mathbb{Z}[x]$ have degrees $k-1$ and $l$ respectively.  Thus by
induction hypothesis,   both
$$P_i(x)=quo(S_{i,  k-1}(x),  l_i(x))\in R$$
and
$$r_i(x)=rem(S_{i,  k-1}(x),  l_i(x))\in R$$
have  a period $b_l^{h}$.   For $i=0,  \cdots,  b_l-1$,   when
$n=b_lm+i$, take $
  P(n) =a_kb_l^{h}m^{h+1}+P_i(m)$ and $
  r(n) =r_i(m)$.
Then it is easy to check that $P(x)=quo(f(x),g(x))$ and
$r(x)=rem(f(x),g(x))$ are elements in $R$ with a period $b_l^{h+1}$,
as desired.

\par \emph{Step 2}: Now we turn to the proof for the uniqueness of $P(x),  r(x)$.  Suppose that there exist $P_1(x),  P_1(x),  r_1(x),  r_2(x)\in R$
such that$$f(x)=P_i(x)g(x)+r_i(x),  \,  0\preccurlyeq  r_i(x)\prec
|g(x)|\,
 \,  \,   i=1,  2$$ Let $T$ be a common
period for $P_1(x),P_2(x),r_1(x)$ and $r_2(x)$.   For $j=0, 1,
\cdots, T-1$), when $n=Tm+j$, we have $P_i(n)=P_{ij}(m)$,
$r_i(n)=r_{ij}(m), \, \, i=1, 2$. Thus
$$(P_{1j}(m)-P_{2j}(m))g(Tm+j)=r_{1j}(m)-r_{2j}(m)$$

By Definition $\ref{2. 0. 2}$,   there exists an integer $C$ such
that for every integer $m>C$,   we have $0\leqslant
r_{ij}(m)<|g(Tm+j)|$. So
$$|r_{1j}(m)-r_{2j}(m)|<|g(Tm+j)|$$
For every $m\in \mathbb{Z}$,   we have $P_{1j}(m)-P_{2j}(m)\in
\mathbb{Z}$.  Therefore,   we have $r_{1j}(x)=r_{2j}(x)$ and
$P_{1j}(x)=P_{2j}(x)$ for every $j=0,  1,  \cdots,  T-1$. \qed
\end{pf}

\begin{rem}\label{2. 0. 01}This division above, which will be called generalized Euclidean division over  $\mathbb{Z}[x]$, almost coincides with
division in $\mathbb{Z}$ pointwisely in the following sense.  Let
$f(x),  g(x)\in \mathbb{Z}[x]$ and
$$f(x)=P_1(x)g(x)+r_1(x),  \,   where \,  0 \preccurlyeq r_1(x)\prec |g(x)|$$By
Definition $\ref{2. 0. 2}$,   the inequality $0\leqslant
r_1(x)<|g(x)|$ provides us  an integer $C_1$,   such that for all $n
>C_1$,
$$\big[\frac{f(n)}{g(n)}\big]=P_1(n)\,  ,  \,  \big\{\frac{f(n)}{g(n)}\big\}g(n)=r_1(n)$$
Put $f_1(x)=f(-x)$ and $g_1(x)=g(-x)$.  Since
$$f_1(x)=P_2(x)g_1(x)+r_2(x),  \,  0\preccurlyeq r_2(x)\prec |g_1(x)|$$
there exists an integer $C_2$,   such that for all $n>C_2$,
$0\leqslant r_2(n)<|g_1(n)|$.  Note
that$$f(-x)=g(-x)P_2(x)+r_2(x)$$Put $y=-x$,   we have
$$f(y)=g(y)P_2(-y)+r_2(-y)$$For every integer $y<-C_2$,   we get
$0\leqslant r_2(-y)<|g(y)|$.  Thus for every integer $n<-C_2$,   we
have
$$\big[\frac{f(n)}{g(n)}\big]=P_2(-n),  \,  \,  \,  \,  \big\{\frac{f(n)}{g(n)}\big\}g(n)=r_2(-n)$$\par
In the special case where $r_1(x)=0$,   we also have $r_2(x)=0$.
Then for {\bf every } $n\in \mathbb{Z}$,
$$r_1(n)=r_2(n)=\big\{\frac{f(n)}{g(n)}\big\}g(n)=0$$
\end{rem}

\begin{exmp}\label{eg1}The following is  an example to
illustrate the relation between generalized  division on
$\mathbb{Z}[x]$ and on $\mathbb{Z}$.
 When $n>1$,
$$\big[\frac{n^{2}}{2n+1}\big]=\left\{
                             \begin{array}{ll}
                               m-1,   & \hbox{$n=2m$;} \\
                               m-1,   & \hbox{$n=2m-1$. }
                             \end{array}
                           \right. $$
$$\big\{\frac{n^{2}}{2n+1}\big\}(2n+1)=\left\{
                             \begin{array}{ll}
                               3m+1,   & \hbox{$n=2m$;} \\
                               m,   & \hbox{$n=2m-1$. }
                             \end{array}
                           \right.
$$
When $n<-1$,  $$\big[\frac{n^{2}}{2n+1}\big]=\left\{
                             \begin{array}{ll}
                               m,   & \hbox{$n=2m$;} \\
                               m,   & \hbox{$n=2m-1$. }
                             \end{array}
                           \right. $$
$$\big\{\frac{n^{2}}{2n+1}\big\}(2n+1)=\left\{
                             \begin{array}{ll}
                               -m,   & \hbox{$n=2m$;} \\
                               -3m+1,   & \hbox{$n=2m-1$. }
                             \end{array}
                           \right.
$$
\end{exmp}

So far, we have defined generalized Euclidean Division on
$\mathbb{Z}[x]$. Now we consider the  case when $f(x),g(x)\in R$.
Suppose that $T_0$ is the least common period of $f(x)$, $g(x)$,
such that $f(x)=(T_0, \{f_i(x)\}_{i=0}^{T_0-1})$ and $g(x)=(T_0,
\{g_i(x)\}_{i=0}^{T_0-1})$.  Based on generalized division on
$\mathbb{Z}[x]$ (in Theorem $\ref{2. 0. 3}$),   we can define
$quo(f(x),  g(x))$ and $rem(f(x), g(x))$ as follows (denoted by
$P(x)$ and $r(x)$ respectively): when $n=Tm+i$
$$P(n)=\left\{
         \begin{array}{ll}
           quo(f_i(m),  g_i(m)),   & \hbox{ if $g_i(m)\neq 0$} \\
           0,   & if \hbox{$g_i(m)=0$}
         \end{array}
       \right.
r(n)=\left\{
         \begin{array}{ll}
           rem(f_i(m),  g_i(m)),   & \hbox{if $g_i(m)\neq 0$} \\
           f_i(m),   &if  \hbox{$g_i(m)=0$}
         \end{array}
       \right. $$

Then it is easy to check that $P(x),   r(x)\in R$ and
\begin{equation}\label{fg}
f(x)=R(x)g(x)+r(x)
\end{equation}
 This will be called  the generalized Euclidean algorithm in the
ring of integer-valued quasi-polynomials.

\section{Generalized
GCD through Generalized Euclidean Division}

When studying Euclidean division in $\mathbb{Z}$,   we can develop
GCD theory and related algorithm (see \cite{Bhu}).  Its counterpart
in $R$ is generalized GCD through generalized Euclidean division.

\begin{defn}[Divisor]\label{2. 0. 5}Suppose that $f(x),  g(x)\in R$ and for every $n\in \mathbb{Z}$,   $g(n)\neq 0$.  Then by Remark
$\ref{2. 0. 01}$ and generalized Euclidean division,   the following
two statements are equivalent:
\begin{enumerate}
  \item[(1)]rem($f(x),  g(x)$)=0;
  \item[(2)] for every $x\in \mathbb{Z}$,   $g(x)$ is a divisor of $f(x)$.
\end{enumerate}
If  the two conditions are satisfied,   we shall call $g(x)$ a
divisor of $f(x)$ and write $g(x)\mid f(x)$.
\end{defn}

By the equivalence of the above two statements,   similar to the
situation in $\mathbb{Z}$,   we have the following proposition.
 \begin{prop}\label{2. 0. 6}Let
$g(x),   f(x)\in $R.  If $f(x)\mid g(x)$ and $g(x)\mid f(x)$,   we
have $f(x)=\varepsilon g(x)$,   where $\varepsilon $ is an
invertible element in $R$.
\end{prop}

\begin{defn}\label{2. 0. 006}[GCD]Let $f_1(x),  f_2(x),  \cdots,  f_s(x),  d(x)\in
R$.
\begin{enumerate}
  \item[(1)] We call $d(x)$  a
common divisor of $f_1(x),  f_2(x),  \cdots,  f_s(x)$,   if we have
$d(x)\mid f_k(x)$ for every $k=1,  2,  \cdots,  s$.
  \item[(2)]We call
$d(x)$  a greatest common divisor  of $f_1(x),  f_2(x),  \cdots,
f_s(x)$ if $d(x)$ is a common divisor  of $f_1(x),  f_2(x),  \cdots,
f_s(x)$ and  for any common divisor $p(x)\in R$ of $f_1(x), f_2(x),
\cdots, f_s(x)$,  we have $p(x)\mid d(x)$.
\end{enumerate}
\end{defn}
\begin{rem}\label{2. 0. 014} Suppose that both $d_1(x)$ and   $ d_2(x)$  are
 greatest common divisors of $f_1(x),  f_2(x), \cdots,  f_s(x)$.
  Then we have $d_1(x)\mid d_2(x)$ and $d_2(x)\mid d_1(x)$. Thus, by
Proposition $\ref{2. 0. 6}$,   we have $d_1(x)=\varepsilon d_2(x)$,
where $\varepsilon $ is an invertible element in $R$.  So we have a
unique GCD $d(x)\in R$ for $f_1(x), f_2(x),  \cdots,  f_s(x)$ such
that $d(x)\succ 0$ or $d(x)=0$ and write it as $ggcd(f_1(x), f_2(x),
\cdots, f_s(x))$.
\end{rem}
\begin{lem}\label{2. 0. 007}Let $f_1(x),  f_2(x),  \cdots,  f_s(x)\in R$.  Then
$$ggcd(f_1(x),  f_2(x),  \cdots,  f_s(x))=ggcd(r(x),  f_2(x),  \cdots,  f_s(x))$$where
$r(x)=rem(f_1(x),  f_2(x))$.
\end{lem}
\begin{pf} By ($\ref{fg}$), $f_1(x)=f_2(x)P_1(x)+r(x) $, we have the set of common
divisors of $f_1(x),  f_2(x),  \cdots,  f_s(x)$ and that of $r(x),
f_2(x),  \cdots,  f_s(x)$ are the same.  The proof is completed.
\qed \end{pf}

\begin{lem}\label{2. 0. 008}Let $f_0(x),  g_0(x)\in R$.  By  generalized Euclidean division,   define
$f_k(x),  g_k(x)(k\in \mathbb{N}-\{0\})$ recursively as follows:
$$f_{k}(x)=g_{k-1}(x),  \,  \,  \,  \,  g_{k}(x)=rem(f_{k-1}(x),  g_{k-1}(x))$$
Then there exists $k_0\in \mathbb{N}-\{0\}$ such that:
$rem(f_{k_0}(x), g_{k_0}(x))=0$.
\end{lem}

\begin{pf}
When $k\geqslant 1$,   if $g_k(x)=rem(f_{k-1}(x),  g_{k-1}(x))\neq
0$.  Let $(T,  \{g_{ki}(x)\}_{i=0}^{T-1})$ be a representation of
$g_{k}(x)$. Define
\begin{eqnarray*}
    g_{ki}(x)&=& a_lx^{l}+a_{l-1}x^{l-1}+\cdots,  +a_1x+a_0 \\
    g_{k-1,  i}(x)&=& b_sx^{s}+b_{s-1}x^{s-1}+\cdots,  +b_1x+b_0=g_{k-1}(Tx+i)
    \end{eqnarray*}
    where $a_l,  b_s\neq 0$.
Since $0\preccurlyeq   g_{ki}(x)\prec |g_{k-1,  i}(x)|$, there are
four cases:
\begin{enumerate}
                   \item[(1)]  $l=s$ and $a_l<|b_s|$;
                   \item[(2)] $l=s$ and $a_l=|b_s|$;
                   \item[(3)] $l<s$;
                   \item[(4)] $g_{ki}(x)=0$.
                 \end{enumerate}
By the generalized Euclidean algorithm in Section 3 and the
remainder condition $0 \preccurlyeq rem(f(x),  g(x))\prec|g(x)|$. We
can reduce case (1) to case (2),   (3) or (4),  reduce (2) to (3) or
(4) and reduce (3) to (4).  For case (4), however, by generalized
Euclidean division,   if $g_{ki}(x)=0$, we have, for every
$t\geqslant k$, $g_{ti}(x)=0$.  Therefore, for all the four cases,
we can find $k_0\in \mathbb{N}$ such that: $rem(f_{k_0}(x),
g_{k_0}(x))=0$.\qed
\end{pf}

\begin{thm}\label{2. 0. 0099}Let $f_1(x),  f_2(x),  \cdots,  f_s(x)\in R$.
Then there exist $d(x),  u_i(x)\in R$ ($ i=1,  2,  \cdots,  s$),
such that $d(x)=ggcd(f_1(x),  f_2(x),  \cdots,  f_s(x))$ and
$$f_1(x)u_1(x)+f_2(x)u_2(x)+\cdots+f_s(x)u_s(x)=d(x)$$
\end{thm}

\begin{pf} To prove this statement,   we apply induction on $s$.  It
is trivial when $s=1$.  Suppose the statement true for $s=l$.  Then
for $s=l+1$,   by induction hypothesis,   there exist $d_0(x), d(x),
u_0(x),  u_1(x),  \cdots,  u_{l+1}(x),  v_1(x),  v_2(x)\in R$ such
that$$d_0(x)=ggcd(f_1(x), f_2(x))=v_1(x)f_1(x)+v_2(x)f_2(x)$$and
$$d(x)=ggcd(d_0(x),  f_3(x),  \cdots,  f_{l+1}(x))=u_0(x)d_0(x)+u_3(x)f_3(x)+\cdots+u_{l+1}(x)f_{l+1}(x)$$
By Lemma $\ref{2. 0. 007}$ and Lemma $\ref{2. 0. 008}$,   we get
$$ggcd(f_1(x),  f_2(x),  f_3(x),  \cdots,  f_{l+1}(x))=ggcd(0,  d_0(x),  f_3(x),  \cdots,  f_{l+1}(x))=d(x)$$
Let $u_1(x)=v_1(x)u_0(x)$ and $ u_2(x)=v_2(x)u_0(x)$.  Then we have
$u_1(x),  u_2(x),  \cdots,  u_{l+1}(x)\in R$ such that
$$d(x)=u_1(x)f_1(x)+u_2(x)f_2(x)+u_3(x)f_3(x)+\cdots+u_{l+1}(x)f_{l+1}(x)$$
\qed \end{pf}
\begin{cor}\label{cor1}Every finitely generated ideal in $R$ is a principal ideal.
\end{cor}
\begin{pf} The result follows readily from Theorem \ref{2. 0. 0099}.
\qed \end{pf}

\section{ GGCD through Pointwise GCD}
Now we are going to study generalized GCD by considering pointwisely
defined GCD in $\mathbb{Z}$ of $f_1(x),  f_2(x),  \cdots,
f_s(x)$($x\in \mathbb{Z}$),   which will give us a more efficient
algorithm for generalized GCD of elements in $R$ than successive
divisions in practice.  \par

Let $f_k(x)\in \mathbb{Z}[x]$,   $k=1,  2,  \cdots,  s$.  For any
$n\in \mathbb{Z}$,   denote the  GCD of $f_1(n),  f_2(n),  \cdots,
f_s(n)$ by $gcd_{\mathbb{Z}}$ $(f_1(n),  f_2(n),  \cdots,  f_s(n))$.
Denote the greatest common factor as polynomials over $\mathbb{Z}$
by $gcd_{\mathbb{Z}[x]}(f_1(x),  f_2(x),  \cdots,  f_s(x))$.

\begin{lem}\label{2. 0. 012}Let $f_i(x)\in \mathbb{Z}[x]$ and
$gcd_{\mathbb{Z}[x]}(f_1(x),  f_2(x),  \cdots,  f_s(x))=1$.  Then
there exists $a_0\in \mathbb{N}$,   such that for all $n \equiv i\,
(mod\, a_0)$,
$$gcd_{\mathbb{Z}}(f_1(n),  f_2(n),  \cdots,  f_s(n))=gcd_{\mathbb{Z}}(f_1(i),  f_2(i),
\cdots,  f_s(i))$$ where $i=0,  1,  \cdots,  a_0-1$.
\end{lem}
\begin{pf} First,
               there exist  $u_k(x)\in \mathbb{Q}[x](k=1,  2,  \cdots,  s)$
               such that $\sum_{k=1}^{s}f_k(x)u_k(x)=1$.  Choose $a_0\in \mathbb{N}$ large enough such that for all $k=1,  2,  \cdots,  s$,
               $w_{k}(x)=a_0u_k(x)\in \mathbb{Z}[x]$.  Since
               $\sum_{k=1}^{s}f_k(x)w_{k}(x)=a_0$,
               we have,   for every $ n\in \mathbb{Z}$,
               \begin{eqnarray}
                gcd_{\mathbb{Z}}(f_1(n),  f_2(n),  \cdots,  f_s(n))\mid a_0\label{a0}
                \end{eqnarray}
                 Then we will show that
                for every $ n\in \mathbb{Z}$,   if $n \equiv i(mod a_0),  i=0,  1,  \cdots,  a_0-1$,   we will have
              $$gcd_{\mathbb{Z}}(f_1(n),  f_2(n),  \cdots,  f_s(n))=gcd_{\mathbb{Z}}(f_1(i),  f_2(i),  \cdots,  f_s(i))$$
                For $i=0,  1,  \cdots,  a_0-1$ and any $m\in \mathbb{Z}$, let
                 \begin{eqnarray}
                d_i&=& gcd_{\mathbb{Z}}(f_1(i),  f_2(i),  \cdots,  f_s(i)) \label{di}\\
                d_{im} &=& gcd_{\mathbb{Z}}(f_1(a_0m+i),  f_2(a_0m+i),  \cdots,  f_s(a_0m+i))\label{dim}
                \end{eqnarray}

                For each $k=1,  2,  \cdots,  s$,   expand $f_k(a_0m+i)$ by binomial expansion theorem.  By
                Lemma $\ref{2. 0. 000}$,
                we get polynomials $h_{ik}(x)\in \mathbb{Z}[x]$ such that
                \begin{eqnarray}
                f_k(n)=f_k(a_0m+i)=a_0h_{ik}(m)+f_k(i)\label{a1}
                \end{eqnarray}
                By ($\ref{di}$),   we have $d_i\mid
                f_k(i)$ for every $k=1,  2,  \cdots,  s$. By ($\ref{a0}$),    we have $d_i\mid a_0$.
                Therefore,   by ($\ref{a1}$),   for all  $m\in \mathbb{Z}$,   we have $d_i\mid
                f_k(a_0m+i)$.  Then from ($\ref{dim}$),   we know $d_i\mid
                d_{im}$.  Similarly,   $d_{im}\mid
                d_i$.  The proof is completed.
\qed\end{pf}
\begin{cor}\label{2. 0. 013}Let $f_1(x),  f_2(x),  \cdots,  f_s(x)\in R$ and for every $n\in \mathbb{Z},  \sum_{i=1}^{s}f_i(n)^{2}\neq
0$. Define a function $g(x)$ as follows:
\begin{eqnarray*}
  g: \mathbb{Z}& \rightarrow & \mathbb{Z} \\
  n &\longmapsto & gcd_{\mathbb{Z}}(f_1(n),  f_2(n),  \cdots,  f_s(n))
\end{eqnarray*}Then $g(x)\in R$.
\end{cor}
\begin{pf} Let $T$ be a common period of
$f_1(x),  f_2(x),  \cdots,  f_s(x)$.  Then there exist $f_{ij}(x)\in
\mathbb{Z}[x]$,   $i=1,  2,  \cdots,  s$,   $j=0,  1,  \cdots,  T-1$
such that if $n=Tm+j$,   $f_i(n)=f_{ij}(m)$.  Put
$$d_{j_0}(x)=gcd_{\mathbb{Z}[x]}(f_{1j}(x),  f_{2j}(x),  \cdots,  f_{sj}(x))$$
By assumption,   we have $d_{j_0}(x)\neq 0$.  Let
$\overline{f_{ij}(x)}=\frac{f_{ij}(x)}{d_{j_0}(x)}$.  Since
$$gcd_{\mathbb{Z}[x]}(\overline{f_{1j}(x)},  \overline{f_{2j}(x)},  \cdots,  \overline{f_{sj}(x)})=1$$by
Lemma $\ref{2. 0. 012}$,
$$d_j(x)=gcd_{\mathbb{Z}}(f_{1j}(x),  f_{2j}(x),  \cdots,  f_{sj}(x))=
d_{j_0}(x)gcd_{\mathbb{Z}}(\overline{f_{1j}(x)},
\overline{f_{2j}(x)},  \cdots,  \overline{f_{sj}(x)})\in R$$Then for
$j=0,  1,  \cdots,  T-1$,  when $n=Tm+j$,   we have
$$g(n)=gcd_{\mathbb{Z}}(f_1(n),  f_2(n),  \cdots,  f_s(n))=d_j(m)$$Thus,
it is easy to check that $g(x)\in R$. \qed \end{pf}
\begin{thm}\label{2. 0. 015}Let $f_1(x),  f_2(x),  \cdots,  f_s(x)\in R$.  Then we have (see Definition $\ref{2. 0. 006}$ and Corollary
$\ref{2. 0. 013}$)$$gcd_{\mathbb{Z}}(f_1(x),  f_2(x),  \cdots,
f_s(x))=ggcd(f_1(x),  f_2(x),  \cdots,  f_s(x))$$
\end{thm}
\begin{pf} Let
$g(x)=gcd_{\mathbb{Z}}(f_1(x),  f_2(x),  \cdots,  f_s(x))$. By
Definition $\ref{2. 0. 5}$,   $g(x)$ is a common divisor of $f_1(x),
f_2(x),  \cdots,  f_s(x)$.  For any common divisor $p(x)\in R$ of
$f_1(x),  f_2(x),  \cdots,  f_s(x)$ we have $p(x)\mid g(x)$. Then
from Definition  $\ref{2. 0. 014}$(ggcd), we should only check that
$g(x)$ is nonnegative, which is obvious. The proof is completed.
\qed
\end{pf}
\begin{exmp}\label{eg2}
From the theory above in this section, we can easily compute ggcd.
For example, $ggcd(x^{3}+2, 3x^{2}-3x, 7x)=h(x)$ where $ h:
\mathbb{Z} \rightarrow  \mathbb{Z} $ and
$$h(n)=\left\{
                              \begin{array}{ll}
                                2,   & \hbox{$n=2m$} \\
                                1,   & \hbox{$n=2m+1$}
                              \end{array}
                            \right.
$$
\end{exmp}

\section{Applications}
\subsection{Application of  Euclidean Division}
Now we will apply generalized successive Euclidean division to
expand  rational numbers with an integer parameter into finite
simple continued fractions(see \cite{Number}).

\begin{thm}\label{finite}Let $f(x),  g(x)\in \mathbb{Z}[x]$ and for every $n\in
\mathbb{Z}$,   $g(n)\neq 0$.    Then there exists $L\in \mathbb{N}$
such that for every $n\in \mathbb{Z}$,   the number of terms  in
expansion $h(n)=\frac{f(n)}{g(n)}$ as finite simple continued
fraction  is no greater than $ L$.
\end{thm}
\begin{pf}By generalized Euclidean division, there exist nonnegative integers $C_1, C_2$ such that $h(n)$ has
uniform expansion formulas ( $n$ being divided into finite cases) as
finite simple continued fractions for all but finitely many values
of $n$, i.e., those integers in the interval $[-C_2,C_1]$. It
follows readily that the numbers of terms for $h(n)$ are
bounded.\qed
\end{pf}

\begin{exmp}Here is an example to illustrate Theorem $\ref{finite}$. \\When $n>4$
$$\frac{n^2}{2n+1}=\left\{
                     \begin{array}{ll}
                       $[m-1;1,  3,  m]$,   & \hbox{$n=2m$} \\
                       $[m-1;3,  1,  m-1]$,   & \hbox{$n=2m-1$}
                     \end{array}
                   \right.
$$
When $n<-4$
$$\frac{n^2}{2n+1}=\left\{
                     \begin{array}{ll}
                       $[m-1;1,2,1,-m-1]$,   & \hbox{$n=2m$} \\
                       $[m-1;4,-m]$,   & \hbox{$n=2m-1$}
                     \end{array}
                   \right.
$$
 \end{exmp}

Similar with finite simple continued fraction,   we can apply
generalized Euclidean division on $\mathbb{Z}[x]$ to some other
problems in elementary number theory,   such as computing Jacobi
Symbol with a parameter in $\mathbb{Z}$.

\begin{exmp}We have\begin{displaymath}
\bigg(\displaystyle\frac{4n^{2}+1}{2n+1}\bigg)=\bigg(\displaystyle\frac{2}{2n+1}\bigg)=\left\{
                                                           \begin{array}{ll}
                                                             1,   & \hbox{$n=4m$} \\
                                                             -1,   & \hbox{$n=4m+1$} \\
                                                             -1,   & \hbox{$n=4m+2$} \\
                                                             1,   & \hbox{$n=4m+3$}
                                                           \end{array}
                                                         \right.
\end{displaymath}
\end{exmp}
\subsection{Application of  generalized GCD}

Based on generalized GCD theory,   we have the following
applications to ideals and matrices in $\mathbb{Z}[x]$.
\begin{lem}\label{2. 0. 4}Let $I_1,  I_2$ be ideals in $\mathbb{Z}[x]$,   then
for every $n\in \mathbb{Z}$,   $I_1(n)=I_2(n)$ as ideals in
$\mathbb{Z}$ if and only if $ggcd(I_1)=ggcd(I_2)$.
\end{lem}
\begin{rem}\label{rem1}Consider the relation between $I_1=I_2$ and
$ggcd(I_1)=ggcd(I_2)$.  It is clear that $I_1=I_2$ implies
$ggcd(I_1)=ggcd(I_2)$.  However,   the converse is not true.  For
example $$ggcd(2,  x+1)=ggcd(4,  x^{2}+1)= \left\{
  \begin{array}{ll}
    2,   & \hbox{$x=2m+1$} \\
    1,   & \hbox{$x=2m$}
  \end{array}
\right.
$$but $(2,  x+1)\neq (4,  x^{2}+1)$.
\end{rem}

\begin{defn}\label{3. 0. 5}Let $A(x)=(a_{ij}(x))_{m\times n},  \,  \,  (m\leqslant n)$,   where $a_{ij}(x)\in
\mathbb{Z}(x)$.  For $k=1,  2,  \cdots,  m$, let $F_k(A(x))$ be the
ideal generated by all the minors of $A(x)$ with order $k$.
\begin{enumerate}
  \item[(1)] Define generalized determinate factors $\{D_k(A(x))\}_{k=1}^{m}$ of
$A(x)$ by $\{ggcd(F_k(A(x)))\}_{k=1}^{m}$;
  \item[(2)]Define generalized invariant factors $\{d_k(A(x))\}_{k=1}^{m}$ by  $d_1(A(x))=D_1(A(x))$ and
$d_k(A(x))=\{quo(D_k(A(x)),  D_{k-1}(A(x)))\}(2\leq k \leq m)$ ;
  \item[(3)] Define the smith normal form $smith(A(x))$ of $A(x)$ of size $m\times n$ by
$$\left(
  \begin{array}{ccccccc}
    d_1(A(x)) & 0 & 0 & 0 & 0&\cdots & 0 \\
    0 & d_2(A(x)) & 0 & 0 & 0&\cdots& 0 \\
    0 & 0 & \cdots & 0 & 0 &\cdots & 0 \\
    0 & 0 & 0 & d_m(A(x)) & 0&\cdots & 0 \\
  \end{array}
\right)_{ m \times n}$$
\end{enumerate}

\end{defn}
By the existence and uniqueness of quotient and ggcd,   the above
three definitions are well defined and  unique for $A(x)$.  Based on
the classical smith normal form for integral matrices (see
\cite{Newmann}),   we have the following result.
\begin{thm}\label{3. 0. 6}Let $A(x)=(a_{ij}(x))_{m\times n},  \,  \,  \,  B(x)=(b_{ij}(x))_{m\times n}$,   where $a_{ij}(x),  b_{ij}(x)\in \mathbb{Z}(x)$.  Then
the following are equivalent:\begin{enumerate}
                               \item [(1)] for $k=1,  2,  \cdots,  m$,   for all $n\in
\mathbb{Z}$,   $ A(n),  B(n)$ are equivalent;

                               \item [(2)]for $k=1,  2,  \cdots,  m$,
$D_k(A(x))=D_k(B(x))$;
                               \item [(3)]for $k=1,  2,  \cdots,  m$,
$d_k(A(x))=d_k(B(x))$;
                               \item [(4)]$smith(A(x))=smith(B(x))$.
                             \end{enumerate}
\end{thm}

\begin{rem}\label{3. 0. 7}By Remark $\ref{rem1}$, it is clear that for each $k=1,  2,  \cdots, m$,
$F_k(A(x))=F_k(B(x))$ implies for every $n\in \mathbb{Z}$, $A(n)$
and $B(n)$ are equivalent and the converse is not true. For example,
let
$$A(x)=\left(
                      \begin{array}{ccc}
                        x & x+2 & 0\\
                        x+1 & x+3 & x \\
                        0 & 0 & 2x\\
                      \end{array}
                    \right)= x
                    \left(
                      \begin{array}{ccc}
                        1 & 1 & 0\\
                        1& 1 & 1 \\
                        0 & 0 & 2\\
                      \end{array}
                    \right)+
                    \left(
                      \begin{array}{ccc}
                        0 & 2 & 0\\
                        1& 3 & 0 \\
                        0 & 0 & 0\\
                      \end{array}
                    \right)
                   $$
 and
$$B(x)=\left(
                      \begin{array}{ccc}
                        x & x+2 & 0 \\
                        x+1 & x+3 & 1 \\
                        0& 0& 2x \\
                      \end{array}
                    \right)=
                    x
                    \left(
                      \begin{array}{ccc}
                        1 & 1 & 0\\
                        1& 1 & 0 \\
                        0 & 0 & 2\\
                      \end{array}
                    \right)+
                    \left(
                      \begin{array}{ccc}
                        0 & 2 & 0\\
                        1& 3 & 1 \\
                        0 & 0 & 0\\
                      \end{array}
                    \right)
$$
Then
$$D_1(A(x))=D_1(B(x))=1$$
\begin{equation*}
D_2(A(x))=D_2(B(x))= \left\{ \begin{aligned}
         1&\,  \,  \,  \,  x\equiv 0(mod 2)\\
         2&\,  \,  \,  \,  x\equiv 1(mod 2)
        \end{aligned} \right.
\end{equation*}
$$D_3(A(x))=D_3(B(x))=4x$$
 According to Theorem $\ref{3. 0. 6}$,   for all $n\in
\mathbb{Z}$,   $A(n)$ and $B(n)$ are equivalent.  But
$$F_2(A(x))=<2,  x^{2}>\neq <2,  x>=F_2(B(x))$$
\end{rem}


\begin{thebibliography}{00}





\bibitem{Atiyah} M. F.  Atiyah.   I. G. Macdonald,   
{Introduction to Commutative Algebra},   {Addison-Wesley Publishing Company},   1969.
\bibitem{stanley} R. P.  Stanley.   
{Enumerative Combinatorics},   Vol.  1.   Cambridge University
Press, 1996.





\bibitem{Number}  K.  H. Rosen. 
{Elementary Number Theory and Its Applications},   Fifth Edition,  
Addison-Wesley Publishing Company,   2004.


\bibitem{Bhu} B.  Mishra.  
{Algorithmic Algebra},  
  Springer-Veriag,   
2001. 

\bibitem{Ring}I. M.  Isaacs.   {Algebra: A Graduate Course},   Wadsworth Inc.  1994.


 \bibitem{Newmann}M. Newmann.    {Integral Matrices},    New York: Academic Press, 1972.
\end{thebibliography}
\end{document}